\documentclass[11pt]{amsart}
\newtheorem{theorem}{Theorem}[section]
\newtheorem{corollary}[theorem]{Corollary}
\newtheorem{lemma}[theorem]{Lemma}

\theoremstyle{definition}
\newtheorem{definition}[theorem]{Definition}
\theoremstyle{remark}
\newtheorem{remark}{Remark}
\numberwithin{equation}{section}

\newcommand{\real}{\mathbb R}

\begin{document}
\title[Nonexistence of nonconstant global minimizers with limit at $\infty$]
{Nonexistence of nonconstant global minimizers with limit at $\infty$ of semilinear elliptic equations in all of $\real^N$}
\author{Salvador Villegas}
\thanks{The author has been supported by the MEC Spanish grants
MTM2006-09282 and MTM2008-00988}
\address{Departamento de An\'{a}lisis
Matem\'{a}tico, Universidad de Granada, 18071 Granada, Spain.}
\email{svillega@ugr.es}


\begin{abstract}

We prove nonexistence of nonconstant global minimizers with limit at infinity of the semilinear elliptic equation $-\Delta u=f(u)$ in the whole $\real^N$, where $f\in C^1(\real)$ is a general nonlinearity and $N\geq 1$ is any dimension. As a corollary of this result, we establish nonexistence of nonconstant bounded radial global minimizers of the previous equation.

\end{abstract}

\maketitle

\section{Introduction and main results}

For $N\geq 1$ and $f\in C^1(\real)$, consider the equation

\begin{equation}\label{mainequation}
-\Delta u=f(u) \ \ \mbox{ in } \real^N.
\end{equation}

We consider classical solutions $u\in C^2(\real^N)$.

Equation (\ref{mainequation}) is the Euler-Lagrange equation associated to the energy functional

$$E_\Omega (u)=\int_\Omega\left( \frac{1}{2}\vert \nabla u\vert^2+G(u)\right) dx\, , \ \ \ \ \ \ \mbox{where }G'=-f$$

\noindent and $\Omega \subset\real^N$ is a bounded domain.

In this note we prove nonexistence of nonconstant global minimizers of (\ref{mainequation}) satisfying $\lim_{\vert x\vert\to \infty}u(x)=u_\infty \in \real$, for every dimension $N\geq 1$. By a global minimizer we mean an absolute minimizer of the energy functional with respect to compactly supported perturbations. In this way, we have

\begin{definition}\label{def}

We say that $u\in C^1(\real^N)$ is a global minimizer of (\ref{mainequation}) if for every smooth bounded domain $\Omega\subset\real^N$ we have

$$E_\Omega (u)\leq E_\Omega (w)\, ,$$

\

\noindent for every function $w\in C^{0,1}(\overline{\Omega})$ such that $w=u$ on $\partial \Omega$.

\end{definition}

It is well known that every global minimizer $u$ of (\ref{mainequation}) is a classical $C^2$ solution of this equation.

Clearly, if we consider  $v\in C^\infty (\real^N)$ with compact support
in a smooth bounded domain $\Omega\subset\real^N$ and take $w=u+t v$ in Definition \ref{def}, we obtain that $t=0$ is an absolute minimizer of the one variable function $t\mapsto E_\Omega (u+tv)$. In particular, the second derivative of this function at $t=0$ is nonnegative, which yields

\begin{equation}\label{stability}
\int_{\real^N} \left( \vert \nabla v\vert^2-f'(u)v^2\right) \,
dx\geq 0.
\end{equation}

A solution $u\in C^2(\real^N)$ of (\ref{mainequation}) satisfying (\ref{stability}) for every $v\in C^\infty (\real^N)$ with compact support in $\real^N$ is called stable. Hence, every global minimizer is a stable solution.

We now consider radial solutions of (\ref{mainequation}). By abuse of notation, we write $u(r)$ instead of $u(x)$, where $r=\vert x\vert$ and $x\in \real^N$. Nonexistence of nonconstant bounded radial stable solutions of (\ref{mainequation}) has been proved by Cabr\'e and Capella \cite{cabrecapella} for every $f\in C^1(\real)$, if $N\leq 8$. The same result holds for $9\leq N\leq 10$, by a result of the author \cite{yo}. Moreover, in \cite{yo} it is established that for every nonconstant bounded radial stable solution of (\ref{mainequation}), there exists $M>0$ such that $\vert u(r)-u_\infty \vert \geq M r^{-N/2+\sqrt{N-1}+2}\, ,\ \forall r\geq 1$, where $u_\infty=\lim_{r\rightarrow \infty}u(r)$ (the existence of this limit is implied by the stability of $u$). On the other hand, for $N\geq 11$, \cite{cabrecapella} constructs a polynomial $f$ which admits a nonconstant bounded radial stable solution of (\ref{mainequation}) (see also \cite{yo} for a large class of nonconstant bounded radial stable solutions of problems of the type (\ref{mainequation}) for $N\geq 11$).

A natural question is whether these examples of nonconstant bounded radial stable solutions of problems of the type (\ref{mainequation}) are also global minimizers (of course, to avoid simple situations, we can always truncate the nonlinearities out of the values of the solutions, and consider bounded nonlinearities). More generally, for which dimensions $N$ (necessarily $N\geq 11$) does problem (\ref{mainequation}) admit a nonconstant bounded radial global minimizer? In this note, we show that the answer to this question is always negative for every dimension $N\geq 1$. In fact, this a consequence of a more general result.

\begin{theorem}\label{noglobal}

Let $f\in C^1(\real)$ and $N\geq 1$. Then equation (\ref{mainequation}) has no nonconstant global minimizers satisfying $\lim_{\vert x\vert\to \infty}u(x)=u_\infty\in\real$.

\end{theorem}

\begin{corollary}\label{noglobalradial}

Let $f\in C^1(\real)$ and $N\geq 1$. Then equation (\ref{mainequation}) has no nonconstant bounded radial global minimizers.

\end{corollary}

The existence of nonconstant  bounded global minimizers with some symmetry properties in certain dimensions is related to the existence of a counter-example of De Giorgi's conjecture \cite{degiorgi}: the level sets of every bounded, monotone in one direction, solution of the Allen-Cahn equation

\begin{equation}\label{allencahn}
-\Delta u=u-u^3 \ \ \mbox{ in } \real^N,
\end{equation}

\noindent must be hyperplanes, at least if $N\leq 8$. The conjecture was proved by Ghoussoub and Gui \cite{GG} in dimension $N=2$, and by Ambrosio and Cabr\'e \cite{ac} in dimension $N=3$. For $4\leq N\leq 8$ and assuming an additional limiting condition on $u$, it has been established by Savin \cite{savin}. Recently, Del Pino, Kowalczyk and Wei \cite{contraejemplo} has found a counter-example of De Giorgi's conjecture for $N\geq 9$, based on a minimal graph $\Gamma$ which is not a hyperplane, found by Bombieri, De Giorgi and Giusti  \cite{bgg} in $\real^N$, $N\geq 9$. Before this recent work, Jerison and Monneau \cite{jm} showed that the existence of a counter-example of the conjecture in $\real^{N+1}$  would be established if one could prove the existence of a global minimizer $v$ of (\ref{allencahn}), with $\vert v\vert<1$, which is even with respect to each coordinate. In \cite{cate}, the instability properties of saddle-shaped solutions of (\ref{allencahn}) in low dimensions are studied, and some computations contained in this work suggest the possibility of finding saddle-shaped global minimizers in higher dimensions, in order to apply Jerison and Monneau's result. On the other hand, Savin \cite{savin} has proved that, for $N\leq 7$, every global minimizer is a function of only one Euclidean variable. Hence, as it is expected, the result of Jerison and Monneau cannot provide a counter-example of the conjecture for $N\leq 8$. Indeed, by a result of Alberti, Ambrosio and Cabr\'e \cite{aac}, the solutions of (\ref{allencahn}) which are monotone functions of only one Euclidean variable, are in fact global minimizers in every dimension.

\section{Proof of the main results}

In this section we prove our main results using
a result of Modica \cite{mod} and a preliminary lemma.
In \cite{mod} Modica proved the following pointwise gradient bound
for global solutions of semilinear elliptic equations.

\begin{theorem}{\bf (Modica \cite{mod})}\label{modica}
Let $G\in C^2(\real)$ be a nonnegative
function and $u$ be a bounded solution of $\Delta
u-G'(u)=0$ in $\real^N$. Then,
\begin{equation}
\frac{\vert\nabla u\vert^2}{2}\leq G(u) \quad {\rm in}\ \real^N.
\end{equation}
In addition, if
$G(u(x_0))=0$ for some  $x_0\in\real^N$, then $u$ is constant.
\end{theorem}

In \cite{mod} this bound was proved under the hypothesis
$u\in C^3(\real^N)$. The result as stated above, which applies to all
solutions ---recall that every solution is $C^{2,\alpha}(\real^N)$
since $G\in C^2(\real)$--- was established in \cite{cafgarseg}.

\begin{lemma}\label{1}

Let $f\in C^1(\real)$ and $N\geq 1$. Let $u$ be a nonconstant global minimizer of (\ref{mainequation}) satisfying $\lim_{\vert x\vert\to \infty}u(x)=u_\infty \in \real$. Then

$$G(s)\geq G(u_\infty) \ \ \forall s\in \real .$$

\end{lemma}

\proof Fix $s\in \real$ and consider $R>1$. Let $w\in C^{0,1}(\overline{B_R})$ the function defined by $w(r)=s$, if $\Vert x\Vert \leq R-1$; and $w(x)=(u(x)-s)(\Vert x\Vert -R+1)+s$, if $R-1<\Vert x\Vert\leq R$. Observe that $w(x)=u(x)$ on $\partial B_R$. Taking into account the boundedness of $u$ and $\vert \nabla u\vert$ in $\real^N$ and evaluating the energy of $w$ in $B_R$, we obtain

$$E_{B_R}(w)=G(s)\frac{\omega_N}{N} (R-1)^N+E_{B_{R}\setminus\overline{B_{R-1}}}(w)\leq
G(s)\frac{\omega_N}{N} (R-1)^N +M R^{N-1}\, ,$$

\noindent where $M>0$ depends only on $u$ and $s$, but not on $R>1$.

By the definition of global minimizer we can assert that

$$\frac{\int_{B_R}G(u)}{R^N}\leq \frac{E_{B_R}(u)}{R^N}\leq
\frac{E_{B_R}(w)}{R^N}\leq
G(s)\frac{\omega_N}{N}\frac{(R-1)^N}{R^N}+\frac{M}{R}$$

Taking limit $R\to\infty$ to the left and right of this expression, we conclude

$$G(u_\infty)\frac{\omega_N}{N}\leq G(s)\frac{\omega_N}{N}.$$

\noindent and the lemma follows. \qed

\

{\bf Proof of Theorem \ref{noglobal}}

Suppose, contrary to our claim, that there exists a nonconstant global minimizer $u$ of (\ref{mainequation}) satisfying $\lim_{\vert x\vert\to \infty}u(x)=u_\infty \in \real$. Without loss of generality we can assume that $\sup u>u_\infty$ (otherwise, replace $u$ by $-u$). Take $x_0\in \real^N$ such that $u(x_0)\geq u(x),\, \forall x\in \real^N$. Choose $G\in C^2(\real)$ such that $G'=-f$ and $G(u_\infty)=0$. From Theorem \ref{modica} and Lemma \ref{1} it follows that $G(u(x_0))>G(u_\infty)=0$. Consider $s_0\in (u_\infty,u(x_0))$ such that $G(s_0)<G(u(x_0))$ and take $s_1\in [s_0,u(x_0)]$ satisfying $G(s_1)=\min_{s\in[s_0,u(x_0)]}G(s)$. It follows immediately that

\begin{equation}\label{yaesta}
s_1\in(u_\infty,u(x_0)),\ \ \ \ \ \ G(s)\geq G(s_1) \ \forall s\in [s_1,u(x_0)].
\end{equation}

Define $\Omega_{s_1}=\left\{ x\in\real^N :u(x)>s_1\right\}$. It is obvious that $\Omega_{s_1}$ is a nonempty bounded open set of $\real^N$. Consider $\Omega$ the connected component of $\Omega_{s_1}$ containing $x_0$. Thus $\Omega$ is a bounded domain satisfying $u\equiv s_1$ on $\partial \Omega$. Therefore, we can consider in Definition \ref{def} the function $w\equiv s_1$ in $\Omega$. Applying (\ref{yaesta}) we obtain

\begin{equation}
E_\Omega (w)=\int_\Omega G(s_1)dx\leq \int_\Omega G(u(x))dx<E_\Omega (u).
\end{equation}

\noindent which contradicts Definition \ref{def}. \qed

\

{\bf Proof of Corollary \ref{noglobalradial}}

We first recall (see \cite[Lem. 2.3]{cabrecapella}) that $u_r$, the usual radial
derivative of $u$, does not vanish in $(0,\infty)$ for every nonconstant radial stable solution of (\ref{mainequation}). Then, since every global minimizer is a stable solution, we deduce the existence of $u_\infty=\lim_{r\rightarrow \infty}u(r)$ for every bounded radial global minimizer $u$ of (\ref{mainequation}). Applying Theorem \ref{noglobal}, the proof is complete.

\

\begin{remark}\label{rem1}

To avoid simple situations, we have not considered constant global minimizers. In fact, it is easily seen that $u=k$ is a global minimizer of (\ref{mainequation}) if, and only if, $G(s)\geq G(k)\, , \forall s\in \real$. The necessary condition follows from Lemma \ref{1}. To see that this condition is also sufficient, consider an smooth bounded domain $\Omega\subset\real^N$ and a function $w\in C^{0,1}(\overline{\Omega})$ (not necessarily satisfying $w=k$ on $\partial \Omega$). Thus $E_\Omega (w)\geq\int_\Omega G(w(x))dx\geq\int_\Omega G(k)dx=E_\Omega (u)$, which is our claim.

\end{remark}

\begin{remark}\label{rem2}

We can consider unbounded radial global minimizers of (\ref{mainequation}). Such solutions exist in every dimension $N\geq 1$. Take, for instance, $f\equiv 1$ and $u_0(x)=-\Vert x\Vert^2/(2N)$. To see that $u_0$ is a global minimizer of $-\Delta u=1$ in $\real^N$, consider a smooth bounded domain $\Omega\subset\real^N$. It is a simple matter to check that the energy functional $E_\Omega (u)=\int_\Omega \left( \vert \nabla u\vert^2 /2-u\right) \,
dx$ has an absolute minimizer in $X=\left\{ u\in H^1(\Omega)\, :u(x)=u_0(x)\, ,\forall x\in \partial\Omega\right\}$. Such an absolute minimizer $w_0$ must satisfy $-\Delta w_0=1$ in $\Omega$; $w_0=u_0$ on $\partial\Omega$. But it is obvious that this problem has an unique solution. Then $w_0=u_0$ in $\Omega$, and consequently $u_0$ is an absolute minimizer of $E_\Omega$ in $X$, which is the desired conclusion.

\end{remark}

{\bf Acknowledgment.} The author would like to thank Xavier
Cabr\'e for suggesting the problem.

\end{document}